\documentclass[12 pt]{amsart}
\usepackage[utf8]{inputenc}
\usepackage{gensymb}
\usepackage{amsmath}
\usepackage{amssymb}
\usepackage{amsthm}
\usepackage{enumitem}
\usepackage{hyperref}
\usepackage[left = 0.8 in, right = 0.8 in, top = 1.25 in, bottom = 1.25 in]{geometry}
\usepackage{graphicx}
\usepackage{textcomp}
\usepackage{blindtext}
\usepackage{hyperref}
\hypersetup{
     colorlinks=true,
     linkcolor= black,
     filecolor= black,
     citecolor = black,      
     urlcolor= black,
     }
\makeatletter
\newcommand\footnoteref[1]{\protected@xdef\@thefnmark{\ref{#1}}\@footnotemark}
\makeatother

\theoremstyle{definition}
\newtheorem{MainTheorem1}{Theorem}
\newtheorem{MainTheorem2}[MainTheorem1]{Theorem}
\newtheorem{theorem}{Theorem}[section]
\newtheorem{lemma}[theorem]{Lemma}

\newtheorem{proposition}[theorem]{Proposition}
\newtheorem{conjecture}[theorem]{Conjecture}

\title{On Sums of Practical Numbers and Polygonal Numbers}
\author{Sai Teja Somu and Duc Van Khanh Tran}
\date{}

\begin{document}


\maketitle

\begin{abstract}
    Practical numbers are positive integers $n$ such that every positive integer less than or equal to $n$ can be written as a sum of distinct positive divisors of $n$. In this paper, we show that all positive integers can be written as a sum of a practical number and a triangular number, resolving a conjecture by Sun.  We also show that all sufficiently large natural numbers can be written as a sum of a practical number and two $s$-gonal numbers.
\end{abstract}

\section{Introduction}

Practical numbers, introduced by Srinivasan in \cite{Srinivasan}, are natural numbers $n$ such that every natural number less than or equal to $n$ can be written as a sum of distinct positive divisors of $n$. In this article, natural numbers refer to positive integers. The sequence of practical numbers can also be found on the Online Encyclopedia of Integer Sequences (OEIS) as entry A005153 \cite{oeis_practical}. Stewart \cite{stewart_practical_charac} and Sierpinski \cite{sierpinski_practical_charac} proved the characterization that a natural number $n\geq 2$ with prime factorization $n = p_1^{a_1}p_2^{a_2}\cdots p_k^{a_k}$, where $p_1 < p_2 < \cdots < p_k$, is practical if and only if $p_1=2$ and
$$p_j \leq \sigma\left(p_1^{a_1}p_2^{a_2}\cdots p_{j-1}^{a_{j-1}}\right)+1$$
for all $2\leq j \leq k$, where $\sigma(\cdot)$ denotes the sum of divisors.

There have been many works on various additive representations of natural numbers involving practical numbers. Melfi \cite{two_conjecture_practical} showed that every even natural number is a sum of two practical numbers. Pomerance and Weingartner \cite{primes_and_practical} proved that every sufficiently large odd number can be written as a sum of a practical number and a prime. Somu et al. \cite{practical_square_quad} proved that all natural numbers congruent to $1$ modulo $8$ are expressible as a sum of a practical number and a square. In this paper, we derive results involving sums of practical and polygonal numbers.

In Section \ref{sec: proof_of_theorem1}, we prove that all natural numbers can be written as a sum of a practical number and a triangular number, resolving the conjecture proposed by Sun in the OEIS entry A208244 \cite{oeis_practical_triangular}. We prove the following theorem.

\begin{MainTheorem1}\label{thrm: prac_plus_triangular}
    Every natural number can be written as a sum of a practical number and a triangular number.
\end{MainTheorem1}

Section \ref{sec: proof_of_theorem2} focuses on additive representations of practical and polygonal numbers more generally. We prove that all sufficiently large natural numbers can be written as a sum of a practical number and two $s$-gonal numbers. We prove the following theorem.

\begin{MainTheorem2}\label{thrm: practical_polygon_sum}
    Let $s$ be a natural number greater than $3$. Then, there exists a natural number $N(s)$ such that all natural numbers greater than $N(s)$ can be written as a sum of a practical number and two $s$-gonal numbers.
\end{MainTheorem2}

Finally, in Section \ref{sec: conjectures}, we propose a couple of conjectures regarding additive representations of natural numbers involving practical numbers. 

\section{Notations}
Besides the usual set-theoretic notations and the notations for inequality, equality, and arithmetic operations, the following notations are used throughout the paper.
{\renewcommand\labelitemi{}
\begin{itemize}
    \item $\mathbb{N}$: This denotes the set of positive integers.
    \item $\mathbb{N}_{0}$: This denotes the set of non-negative integers.
    \item $\mathbb{P}_r$: This denotes the set of practical numbers.
    \item $\sigma(\cdot)$: This denotes the sum of positive divisors of a natural number.
    \item $\equiv$: We say $a \equiv b\, (\text{mod } n)$ if $n$ divides $a-b$.
    \item $(\cdot, \cdot)$: This denotes the greatest common divisor of two integers.
    \item $\lfloor \cdot \rfloor$: This denotes the floor function of a real number, which is the largest integer not exceeding that real number.
    \item $P_s(n)$: For integers $s\geq 3$ and $n\geq 0$, this denotes the $n$-th $s$-gonal number given by
    $$P_s(n) = (s-2) \frac{n(n-1)}{2} + n = \frac{(s-2)n^2 - (s-4)n}{2}.$$
    \item $O(\cdot)$: We say $f(x) = O(g(x))$ if there exists a positive real number $M$ such that $|f(x)| \leq Mg(x)$ for all sufficiently large $x$.
\end{itemize}
Other miscellaneous notations are defined where they are used.
}

\section{Proof of Theorem 1}\label{sec: proof_of_theorem1}

In this section, we prove that all natural numbers can be written as a sum of a practical number and a triangular number, resolving the conjecture proposed by Sun in the OEIS entry A208244 \cite{oeis_practical_triangular}. We require three lemmas to prove Theorem \ref{thrm: prac_plus_triangular}.

\begin{lemma}\label{lem: square_8n+1}
    Let $m$ and $n$ be natural numbers. There exists a natural number $1\leq x \leq 2^m-1$ such that
    $$x^2 \equiv 8n+1\, \left(\text{mod } 2^{m+2}\right).$$
\end{lemma}
\begin{proof}
    See \cite[Lemma 3.2]{practical_square_quad} for proof.
\end{proof}

\begin{lemma}\label{lem: practical_product}
    If $m$ is a practical number and $n$ is a natural number such that $n \leq \sigma(m)+1$, then $mn$ is a practical number.
\end{lemma}
\begin{proof}
    See \cite[Corollary 1]{Margenstern_les_nombres_pratiques} for proof.
\end{proof}

\begin{lemma}\label{lem: triangular_number}
If $x$ is an odd natural number, then $\frac{x^2-1}{8}$ is a triangular number.
\end{lemma}
\begin{proof}
As $x=2k+1$ for some non-negative integer $k$, we have $\frac{x^2-1}{8} = \frac{k(k+1)}{2}$. Hence, $\frac{x^2-1}{8}$ is a triangular number.
\end{proof}

Now we prove Theorem \ref{thrm: prac_plus_triangular}.

\begin{proof}[Proof of Theorem \ref{thrm: prac_plus_triangular}]
    Let $n$ be a natural number and $m = \left\lfloor\log_2{\sqrt{8n+1}}\right\rfloor$. By Lemma \ref{lem: square_8n+1}, there exists a natural number $x$ such that $1\leq x \leq 2^m-1$ and
    $x^2 \equiv 8n+1\, \left(\text{mod } 2^{m+2}\right).$
    Since $x \leq 2^m-1$ and $m = \left\lfloor\log_2{\sqrt{8n+1}}\right\rfloor$, we have
    $x^2 < 2^{2m} \leq 8n+1$.
    
    As $x^2 < 8n+1$ and $x^2 \equiv 8n+1\, \left(\text{mod } 2^{m+2}\right)$,
    $8n+1 - x^2 = 2^{m+2}s$ for some natural number $s$. As $m = \left\lfloor\log_2{\sqrt{8n+1}}\right\rfloor$, we have
    $2^{m+2}s \leq 8n+1 \leq 2^{2m+2},$
    which implies that $s \leq 2^m$. Since $2^{m-1}$ is a practical number and
    $s \leq 2^m = \sigma\left(2^{m-1}\right) + 1,$
    by Lemma \ref{lem: practical_product}, $2^{m-1}s$ is a practical number. Notice that $\frac{x^2-1}{8}$ is a triangular number, as $x$ is an odd natural number because $x^2 \equiv 8n+1\, \left(\text{mod } 2^{m+2}\right)$. Now,  $8n+1 - x^2 = 2^{m+2}s$, or equivalently $n = 2^{m-1}s + \frac{x^2-1}{8}$. Since $2^{m-1}s$ is a practical number and $\frac{x^2-1}{8}$ is a triangular number, $n$ is a sum of a practical number and a triangular number. Therefore, all natural numbers can be written as a sum of a practical number and a triangular number.
\end{proof}

\section{Proof of Theorem 2}\label{sec: proof_of_theorem2}
In this section, we prove some results regarding additive representations of natural numbers involving practical and polygonal numbers more generally. We require five lemmas to prove Theorem \ref{thrm: practical_polygon_sum}.

\begin{lemma}\label{lem: sum_poly_mod_p}
    Let $p$ be an odd prime, $n$ be a natural number, and $s$ be a natural number greater than $3$. There exist natural numbers $x$ and $y$ such that
    $$P_s(x) + P_s(y) \equiv n \, (\text{mod } p).$$
\end{lemma}
\begin{proof}
    If $p\mid (s-2)$, then $x=n$ and $y=p$ satisfy
    $$P_s(x) + P_s(y) \equiv n \, (\text{mod } p).$$
    If $p\nmid (s-2)$, then
    $$\left|\{P_s(i) \, (\text{mod } p): 1 \leq i \leq p\}\right| = \left|\{(n-P_s(j)) \, (\text{mod } p): 1 \leq j \leq p\}\right| = \frac{p+1}{2}.$$
    This implies that
    $$\{P_s(i) \, (\text{mod } p): 1 \leq i \leq p\} \cap \{(n-P_s(j)) \, (\text{mod } p): 1 \leq j \leq p\} \neq \varnothing.$$
    Thus, there exist natural numbers $x$ and $y$ such that
    \begin{align*}
        P_s(x) + P_s(y) &\equiv n \, (\text{mod } p).
    \end{align*}
\end{proof}

\begin{lemma}\label{lem: sum_poly_mod_2}
    Let $n$ be a natural number and $s$ be a natural number greater than $3$. There exist natural numbers $x$ and $y$ such that
    $$P_s(x_1) + P_s(y_1) \equiv n \, (\text{mod } 2)$$
    for all natural numbers $x_1 \equiv x\, (\text{mod } 4)$ and $y_1 \equiv y\, (\text{mod } 4)$.
\end{lemma}
\begin{proof}
    If $n$ is even, then $(x,y)=(4,4)$ satisfies the condition above. If $n$ is odd, then $(x,y)=(4,1)$ satisfies the condition above.
\end{proof}

\begin{lemma}\label{lem: sum_square_mod_p^k}
    Let $p$ be a prime congruent to 1 modulo 4. For any natural numbers $n$ and $k$, there exist natural numbers $x$ and $y$ such that
    $x^2 + y^2 \equiv n\, \left(\text{mod } p^k\right)$
    and $p\nmid y$.
\end{lemma}
\begin{proof}
    We will prove the lemma using mathematical induction on $k$. Let us first prove the lemma for $k=1$. Let $n$ be any natural number. By \cite[Theorem 84]{Hardy_Wright},
    $$|\{i^2\, (\text{mod } p): 1 \leq i \leq p\}| = |\{(n - j^2)\, (\text{mod } p): 1 \leq j \leq p\}| = \frac{p+1}{2},$$
    so
    $$\{i^2\, (\text{mod } p): 1 \leq i \leq p\} \cap \{(n - j^2)\, (\text{mod } p): 1 \leq j \leq p\} \neq \varnothing,$$
    and thus there exist natural numbers $x$ and $y$ such that
    $$x^2 + y^2 \equiv n\, \left(\text{mod } p\right).$$
    If $n \not\equiv 0\, (\text{mod } p)$, $x$ and $y$ cannot both be multiples of $p$. Without loss of generality, we can let $p \nmid y$. Since $p \equiv 1 \, (\text{mod } 4)$, there exists a natural number $a$ such that $a^2 + 1 \equiv 0\, (\text{mod }p)$
    (see \cite[Theorem 86]{Hardy_Wright}). So, $x=a$ and $y=1$ is a solution to
    $x^2 + y^2 \equiv 0\, (\text{mod } p)$ with $p\nmid y$. 
    
    Now suppose that there exist natural numbers $x$ and $y_s$ such that
    $$x^2 + y_s^2 \equiv n\, \left(\text{mod } p^s\right)$$
    and $p \nmid y_s$, where $s \geq 1$. Let $l$ be any natural number satisfying $$l\equiv \left(\frac{n-x^2-y_s^2}{p^s}\right) (2y_s)^{-1} (\text{mod } p),$$
    and let $y_{s+1} = y_s+p^sl$. Now, as $y_{s+1}\equiv y_s \, \left(\text{mod } p^s\right)$ and $p \nmid y_s$, we have $p\nmid y_{s+1}$. As $\frac{x^2+y_s^2-n}{p^s} + 2ly_s\equiv 0\, (\text{mod } p)$, we have
    \begin{align*}
    x^2+y_{s+1}^2 &= x^2+(y_s+p^sl)^2 \\
    &= x^2+y_s^2+2p^sly_s + p^{2s}l^2\\
    &\equiv x^2 + y_s^2 + 2p^sl y_s\, (\text{mod } p^{s+1})\\
    &\equiv n + p^s\left(\frac{x^2+y_s^2-n}{p^s} + 2ly_s\right) (\text{mod } p^{s+1})\\
    & \equiv n\, (\text{mod } p^{s+1}).
    \end{align*}
    Hence, by mathematical induction, for all natural numbers $k$, there exist natural numbers $x$ and $y$ such that
    $$x^2 + y^2 \equiv n\, \left(\text{mod } p^k\right)$$
    and $p\nmid y$.
\end{proof}

\begin{lemma}\label{lem: sum_poly_mod_p^k}
    Let $s$ be a natural number greater than $3$. There exists an odd prime $p$ not dividing $s-2$ such that for all $k,n \in \mathbb{N}$, there exist $x,y \in \mathbb{N}$ such that
    $$P_s(x) + P_s(y) \equiv n\, \left(\text{mod } p^k\right).$$
\end{lemma}
\begin{proof}
    Let $k$ and $n$ be any natural numbers, and let $p$ be a prime congruent to 1 modulo 4 such that $p\nmid (s-2)$. Note that
    $$8(s-2)P_s(x) = (2(s-2)x - (s-4))^2 - (s-4)^2.$$
    Since $p \equiv 1 \, (\text{mod } 4)$, by Lemma \ref{lem: sum_square_mod_p^k}, there exist natural numbers $x_0$ and $y_0$ such that
    $$x_0^2 + y_0^2 \equiv 8(s-2)n + 2(s-4)^2 \left(\text{mod } p^k\right).$$
    Let $x$ and $y$ be natural numbers satisfying the congruences
    $$x \equiv 2^{-1}(s-2)^{-1} (x_0 + s - 4) \, \left(\text{mod } p^k\right)$$
    and
    $$y \equiv 2^{-1}(s-2)^{-1} (y_0 + s - 4) \, \left(\text{mod } p^k\right).$$
    Then,
    \begin{align*}
        8(s-2)P_s(x) + 8(s-2)P_s(y) &\equiv x_0^2+y_0^2 - 2(s-4)^2 \left(\text{mod } p^k\right)\\
        &\equiv 8(s-2)n \left(\text{mod } p^k\right).
    \end{align*}
    Since $\left(p^k, 8(s-2)\right) = 1$, we have $P_s(x) + P_s(y) \equiv n\, \left(\text{mod } p^k\right)$.
\end{proof}

\begin{lemma}\label{lem: constant_upper_bound}
    Let $s$ be a natural number greater than $3$, and let $p_{i(s)}$ be the smallest prime for which Lemma \ref{lem: sum_poly_mod_p^k} holds. There exists a real number $A(s)$ such that for all $x\geq 1$, we have
    $$ \frac{2P_s(2p_{i(s)}x)}{x^2} \leq A(s).$$
\end{lemma}
\begin{proof}
    Since $2P_s(2p_{i(s)}x)$ is a quadratic polynomial, we have
    $\frac{2P_s(2p_{i(s)}x)}{x^2} = O(1).$
    Hence, there exists a real number $A(s)$ such that
    $$\frac{2P_s(2p_{i(s)}x)}{x^2} \leq A(s)$$
    for all real numbers $x \geq 1$.
\end{proof}
Now we are ready to give a proof of Theorem 2.
\begin{proof}[Proof of Theorem 2]
    Let $p_i$ denote the $i$-th prime, and let $p_{i(s)}$ be the smallest prime for which Lemma \ref{lem: sum_poly_mod_p^k} holds. By Lemma \ref{lem: constant_upper_bound}, there exists a real number $A(s)$ such that
    $$\frac{2P_s(2p_{i(s)}x)}{x^2} \leq A(s)$$
    for all real numbers $x \geq 1$. Let $r$ be the smallest natural number such that $r \geq i(s)$ and
    $$\frac{\sigma(p_1p_2\cdots p_r)}{p_1p_2\cdots p_r} \geq A(s).$$
    Such an $r$ is well-defined, as the product
    $\prod_{p \text{ prime}} \left(1+\frac{1}{p}\right)$
    diverges (see \cite[Chapter 7, Theorem 3]{Knopp} and \cite[Theorem 19]{Hardy_Wright}). Let $N(s) = 2P_s(2p_1p_2\cdots p_r)$, and consider any natural number $n$ greater than $N(s)$. Let $k$ be the largest natural number such that
    $$2P_s\left(2p_1p_2\cdots p_{i(s)-1}p_{i(s)}^kp_{i(s)+1} \cdots p_r\right) < n.$$ Let $n_k=2p_1p_2\cdots p_{i(s)-1}p_{i(s)}^kp_{i(s)+1} \cdots p_r$. From the definition of $k$, we have $2P_s(n_k) < n \leq 2P_s(p_{i(s)}n_k).$
    From Lemma \ref{lem: sum_poly_mod_p}, there exists a solution $x\, (\text{mod } p_i)$, $y\, (\text{mod } p_i)$ to the equation $P_s(x)+P_s(y) \equiv n\, (\text{mod } p_i)$ for $2\leq i \leq r$ and $i\neq i(s)$. From Lemma \ref{lem: sum_poly_mod_2}, there exists a solution $x \, (\text{mod } 2p_1)$, $y\, (\text{mod } 2p_1)$ to the equation $P_s(x)+P_s(y)\equiv n\, (\text{mod } p_1)$. From Lemma \ref{lem: sum_poly_mod_p^k}, there exists a solution $x\, \left(\text{mod } p_{i(s)}^k\right)$, $y\, \left(\text{mod } p_{i(s)}^k\right)$ to $P_s(x)+P_s(y)\equiv n\, \left(\text{mod } p_{i(s)}^k\right)$. Hence, by the Chinese Remainder Theorem, there exists a solution $x\, (\text{mod } n_k)$, $y\, (\text{mod } n_k)$ to the equation
    $$P_s(x) + P_s(y) \equiv n \, \left(\text{mod } \frac{n_k}{2}\right).$$
    Thus, there exist natural numbers $x, y \leq n_k$ such that
    $$P_s(x) + P_s(y) \equiv n \, \left(\text{mod } \frac{n_k}{2}\right).$$
    This, together with the fact that $n > 2P_s(n_k)$, implies that
    $\frac{2(n-P_s(x)-P_s(y))}{n_k} \in \mathbb{N}.$
    Note that
    \[
        \frac{2(n-P_s(x)-P_s(y))}{n_k}  \leq \frac{2n}{n_k} 
        \leq \frac{4P_s(p_{i(s)}n_k)}{n_k}.
    \]
    By Lemma \ref{lem: constant_upper_bound},
    $\frac{2P_s(p_{i(s)}n_k)}{\frac{n_k^2}{4}} \leq A(s)$. Hence, \[\frac{2(n-P_s(x)-P_s(y))}{n_k} \leq \frac{4P_s(p_{i(s)}n_k)}{n_k} \leq \frac{A(s)n_k}{2}.\]
    Also, we have
    $$\frac{\sigma(\frac{n_k}{2})}{\frac{n_k}{2}} = \frac{\sigma\left(p_1p_2\cdots p_{i(s)-1} p_{i(s)}^kp_{i(s)+1} \cdots p_r\right)}{p_1p_2\cdots p_{i(s)-1} p_{i(s)}^kp_{i(s)+1} \cdots p_r} \geq \frac{\sigma(p_1p_2\cdots p_r)}{p_1p_2\cdots p_r} \geq A(s),$$
    or equivalently $\frac{A(s)n_k}{2}\leq \sigma\left(\frac{n_k}{2}\right).$
    So,
    \begin{align*}
        \frac{2(n-P_s(x)-P_s(y))}{n_k} &\leq \frac{A(s)n_k}{2} \\
        &\leq \sigma\left(\frac{n_k}{2}\right).
    \end{align*}
    Note that $\frac{n_k}{2} = p_1p_2\cdots p_{i(s)-1} p_{i(s)}^kp_{i(s)+1} \cdots p_r$ is a practical number by the characterization of practical numbers (see \cite[Section 3]{stewart_practical_charac}). Thus, by Lemma \ref{lem: practical_product},
    $$\frac{2(n-P_s(x)-P_s(y))}{n_k}\left(\frac{n_k}{2}\right) = n-P_s(x) - P_s(y)$$
    is a practical number. Therefore, $n$ can be written as a sum of a practical number and two $s$-gonal numbers. 
\end{proof}

In Theorem 1, we have proved that all natural numbers can be written as a sum of a triangular number and a practical number. In Theorem 2, we have proved that for all $s>3$, all sufficiently large natural numbers can be written as a sum of a practical number and two $s$-gonal numbers. Now we show that there are infinitely many $s>3$ for which we cannot write all sufficiently large natural numbers as a sum of a practical number and an $s$-gonal number. We also show that as $s$ tends to infinity, the number of natural numbers that cannot be written as a sum of a practical number and two $s$-gonal numbers tends to infinity. Hence, we cannot drop ``sufficiently large" from the statement of Theorem 2. We will require one lemma to prove these claims.

\begin{lemma}\label{lem: practical_not_divisible_3_4}
    If $q$ is a practical number such that $q$ is not divisible by $3$ and $4$, then $q=1$ or $q=2$.
\end{lemma}
\begin{proof}
    For the sake of contradiction, assume that $q > 2$. Since $q>2$ and $4\nmid q$, $q$ should have at least one odd prime divisor. Let $p$ be the smallest odd prime divisor of $q$. As $3\nmid q$, we have $p\geq 5$. As $p \geq 5 > \sigma(2) + 1$, by the characterization of practical numbers, $q$ is not practical (see \cite[Section 3]{stewart_practical_charac}). This is a contradiction.
\end{proof}

\begin{proposition}\label{prop: prac_plus_s-gonal}{\ \\}
    (a) If $s \equiv 0\, (\text{mod } 12)$ or $s \equiv 4\, (\text{mod } 12)$, then infinitely many natural numbers cannot be written as a sum of a practical number and an $s$-gonal number.\\
    (b) Let $E(s)$ be the number of natural numbers that cannot be written as a sum of a practical number and two $s$-gonal numbers. Then, $$\lim_{s\rightarrow \infty} E(s) = \infty.$$
\end{proposition}
\begin{proof}{\ \\}
    (a) If $s \equiv 0\, (\text{mod } 12)$ or $s \equiv 4\, (\text{mod } 12)$, then
    $$P_s(n) = \frac{s-2}{2}n^2 + \frac{s-4}{2}n = an^2 + bn,$$
    where $a$ is odd, $a$ is not divisible by $3$, and $b$ is even. Hence,
    $$P_s(n) = \frac{1}{a} \left(\left(an + \frac{b}{2}\right)^2 - \frac{b^2}{4}\right),$$
    and thus for all $n\in \mathbb{N}_{0}$,
    $$P_s(n) \not\equiv a^{-1}\left(2 - \frac{b^2}{4}\right)(\text{mod } 3),$$
    and
    $$P_s(n) \not\equiv a^{-1}\left(2 - \frac{b^2}{4}\right)(\text{mod } 4)$$
    since $2$ is a quadratic non-residue modulo $3$ and modulo $4$. Let $r$ be a natural number such that
    $$r \equiv a^{-1}\left(2 - \frac{b^2}{4}\right) (\text{mod } 12),$$
    and let $x$ be any positive real number. Let $k$ be a natural number congruent to $r$ modulo $12$ and not exceeding $x$ such that $k = P_s(m) + q$ for some $m\in \mathbb{N}_{0}$ and some practical number $q$. Since $k \equiv r \, (\text{mod } 12)$, we have $3 \nmid (k - P_s(m))$ and $4 \nmid (k - P_s(m))$. Hence, by Lemma \ref{lem: practical_not_divisible_3_4}, $q=1$ or $q=2$, and thus
    $$k \in \{P_s(n)+1: n\in\mathbb{N}_{0} \text{ and } P_s(n)+1 \leq x\} \cup \{P_s(n)+2: n\in\mathbb{N}_{ 0} \text{ and } P_s(n)+2 \leq x\}.$$
    Since
    $$|\{P_s(n)+1: n\in\mathbb{N}_{0} \text{ and } P_s(n)+1 \leq x\}| + |\{P_s(n)+2: n\in\mathbb{N}_{0} \text{ and } P_s(n)+2 \leq x\}| = O\left(\sqrt{x}\right),$$
    there are at most $O\left(\sqrt{x}\right)$ natural numbers not exceeding $x$ that are congruent to $r$ modulo $12$ and expressible as a sum of a practical number and an $s$-gonal number. Since there are $\frac{x}{12} + O(1)$ natural numbers not exceeding $x$ that are congruent to $r$ modulo $12$, there are at least
    $\frac{x}{12} + O\left(\sqrt{x}\right)$
    natural numbers not exceeding $x$ that are congruent to $r$ modulo $12$ and not expressible as a sum of a practical number and an $s$-gonal number. Therefore, infinitely many natural numbers cannot be written as a sum of a practical number and an $s$-gonal number if $s \equiv 0\, (\text{mod } 12)$ or $s \equiv 4\, (\text{mod } 12)$.\\
    (b) Let us count the number of natural numbers less than $s$ that can be represented as a sum of a practical number and two $s$-gonal numbers. Only $0$ and $1$ are $s$-gonal numbers less than $s$. Thus, if $n<s$ is a sum of a practical number and two $s$-gonal numbers, then \[n\in \{P+i: P\in \mathbb{P}_r, P \leq s, \text{ and } i\in \{0,1,2\}\}.\] From \cite{practical_distribution_divisor}, we have $$\left|\{P+i: P\in \mathbb{P}_r, P \leq s, \text{ and } i\in \{0,1,2\}\}\right| = O\left(\frac{s}{\log s}\right).$$ Hence, at least $s+O\left(\frac{s}{\log s}\right)$ natural numbers less than $s$ cannot be written as a sum of a practical number and two $s$-gonal numbers. Therefore,  $$\lim_{s\rightarrow \infty} E(s) = \infty.$$
\end{proof}

\section{Conjectures on Sums of Practical and Polygonal Numbers}\label{sec: conjectures}

In this section, we propose a few conjectures on some additive representations involving practical and polygonal numbers based on numerical computations. All of the codes for the conjectures below, written in Python, can be found on \href{https://github.com/ducvktran/On-Sums-of-Practical-Numbers-and-Polygonal-Numbers}{GitHub}\footnote{\label{foot: 1}\href{https://github.com/ducvktran/On-Sums-of-Practical-Numbers-and-Polygonal-Numbers}{https://github.com/ducvktran/On-Sums-of-Practical-Numbers-and-Polygonal-Numbers}}.

For each natural number $s>3$, let $n_s$ be the number of natural numbers below $10^8$ that cannot be written as a sum of a practical number and an $s$-gonal number, and let $N_s$ be the largest number below $10^8$ that cannot be written as a sum of a practical number and an $s$-gonal number. We have the following conjecture based on Table \ref{tab: practical_polygonal}.
\begin{table}[h]
    \centering
    \begin{tabular}{|c|c|c|}
        \hline $s$ & $n_s$ & $N_s$ \\
        \hline $4$ & $17929061$ & $99999998$\\
        \hline $5$ & $13$ & $2671$\\
        \hline $6$ & $101$ & $1332329$\\
        \hline $7$ & $73$ & $79445$\\
        \hline $8$ & $414$ & $4005819$\\
        \hline
    \end{tabular}
    \hspace{0.3 cm}
    \begin{tabular}{|c|c|c|}
        \hline $s$ & $n_s$ & $N_s$ \\
        \hline $9$ & $186$ & $325808$\\
        \hline $10$ & $341$ & $13613213$\\
        \hline $11$ & $68$ & $105712$\\
        \hline $12$ & $16663689$ & $99999998$\\
        \hline $13$ & $609$ & $1612172$\\
        \hline
    \end{tabular}
    \hspace{0.3 cm}
    \begin{tabular}{|c|c|c|}
        \hline $s$ & $n_s$ & $N_s$ \\
        \hline $14$ & $79$ & $106878$\\
        \hline $15$ & $767$ & $1486748$\\
        \hline $16$ & $16665797$ & $99999998$\\
        \hline $17$ & $106$ & $9314$\\
        \hline $18$ & $1020$ & $8541224$\\
        \hline
    \end{tabular}
    \vspace{0.2 cm}
    \caption{Data regarding the sum of a practical number and an $s$-gonal number}
    \label{tab: practical_polygonal}
\end{table}

\begin{conjecture}
    If $s>3$, $s \not\equiv 0\, (\text{mod } 12)$, and $s \not\equiv 4\, (\text{mod } 12)$, then all sufficiently large natural numbers can be written as a sum of a practical number and an $s$-gonal number.
\end{conjecture}

For $s \in \{4,5,6,7,8,10\}$, we also computationally verified that all natural numbers below $10^8$ are expressible as a sum of a practical number and two $s$-gonal numbers. We have the following conjecture, which is a stronger version of Theorem \ref{thrm: practical_polygon_sum}, for the cases $s \in \{4,5,6,7,8,10\}$.

\begin{conjecture}\label{conj: practical + 2 polygonal}
    For $s \in \{4,5,6,7,8,10\}$, all natural numbers can be written as a sum of a practical number and two $s$-gonal numbers.
\end{conjecture}

The conjecture above does not hold for other values of $s$. One can verify that $23$ cannot be written as a sum of a practical number and two nonagonal numbers and that $11$ cannot be written as a sum of a practical number and two $s$-gonal numbers if $s\geq 11$.

\section{Disclosure Statement}
There is no potential conflict of interest that could influence this paper.

\section{Data Availability Statement}
The relevant data for this study are available within the article.

\bibliography{bibliography.bib}
\bibliographystyle{abbrv}

\end{document}